\documentstyle[12pt]{amsart}
\newtheorem{introteo}{Theorem}

\newtheorem{Thm}{Theorem}[section]
\newtheorem{Prop}[Thm]{Proposition}
\newtheorem{Lem}[Thm]{Lemma}
\newtheorem{case}{Case}
\newtheorem{Cor}[Thm]{Corollary}
\theoremstyle{definition}
\newtheorem{claim}{Claim}
\newtheorem{Def}{Definition}

\newtheorem{rem}[Thm]{Remark}
\begin{document}
\title[Numerical Criteria]{Numerical Criteria for Very Ampleness of
Divisors
on Projective Bundles over an
 Elliptic Curve}
\author[A. Alzati]{Alberto Alzati }
\address{ Dipartimento di Matematica\\ Universita' degli Studi di Milano\\
Italy}
\email{alzati@@vmimat.mat.unimi.it}
\author[M. Bertolini]{ Marina Bertolini}
\address{ Dipartimento di Matematica\\ Universita' degli Studi di Milano\\
Italy}
\email{bertolin@@vmimat.mat.unimi.it}
\author[G. Besana]{ Gian Mario Besana}
\address {Department of Mathematics\\Eastern Michigan
University\\Ypsilanti MI 48197}
 \email{gbesana@@emunix.emich.edu}
\thanks{The authors acknowledge support from GNSAGA, CNR, Italy}
\maketitle
\begin{abstract}
Let $D$ be a divisor on a projectivized bundle over an elliptic curve.
Numerical
conditions for the very ampleness of $D$ are proved. In some cases a
complete
numerical characterization is found.
\end{abstract}
\section{Introduction}
Ampleness of  divisors on algebraic varieties is a numerical property. On
the
other hand it is in general very difficult to give numerical necessary and
sufficient
conditions for the very ampleness of divisors. In \cite{bu} the author
gives a
sufficient condition for a line bundle associated with a  divisor $D$ to
be normally generated on $X
=\Bbb{P}(E)$ where $E$ is a vector bundle over a smooth curve $C.$
A line bundle which is ample and normally generated is automatically very
ample.
Therefore the condition   found in \cite{bu}, together with Miyaoka's well
known
ampleness criterion, give  a sufficient condition for the very ampleness
of $D$ on $X.$
This work is devoted to  the study of numerical criteria for very
ampleness of
divisors $D$ which do not satisfy the above criterion, in the case of $C$
elliptic.
With this assumption Biancofiore and Livorni \cite{bi-li3} (see also
\cite[Prop.8.5.8]{BESO} for a generalization)
 gave a necessary and
sufficient condition when $E$ is indecomposable,
 rk$E$ = 2 and deg$E$ = 1. Gushel \cite{gu}  also gave a
complete characterization of the very ampleness of $D$  assuming that $E$
is
indecomposable and
$|D|$ embeds $X$ as  a scroll. This work deals with the general situation
and addresses
the cases still open.

The main technique used here is a very classical one. A suitable divisor
$A$ on $X$ is
chosen such that there exists a smooth $S \in |A|$ containing  every pair
of points,
possibly infinitely near. Appropriate vanishing conditions are established
to assure
that the natural restriction map $H^0(X, \cal{O}_X(D)) \to H^0(S,
\cal{O}_X(D)_{|S})$ is surjective. In this
 way we get that a divisor $D$ of $X$ is very ample if and only if
$D_{|S}$ is very ample.
In this context $S$ is chosen as $S = \Bbb{P}(E')$ where $E'$ is a
quotient of
$E,$ thus with rank smaller than rank $E.$ Therefore an inductive process
on the
rank can be set up. This process is not always easy to
carry on. For example if $E$  is
assumed to be
 indecomposable there is no guarantee that $E'$ will still be
indecomposable. Since
ampleness is inherited by quotients, we will require at some stage that
$E$ be ample.
The paper is organized as follows. Section 2 contains notation, known and
preliminary results used in the sequel.  In section 3 the case of rank $E
= 2$ is fully
treated. We recover Biancofiore and Livorni's results and deal with the
case of
$E$  decomposable. Section 4 deals with the case of rank $E = 3$ while
section 5
contains the study of case  rank $E \ge 4.$
In particular in the case of rank $E = 3$ we get the following result (see
section \ref{notation} for notation):
\begin{introteo}
Let $E$ be a rank $3$ vector bundle on an elliptic curve $C$ and let
$D\equiv aT+bf$
be a line bundle on  $ X =\Bbb{P}(E).$
\begin{itemize}
\item[(a)] If $E$ is indecomposable then
	\begin{itemize}
	\item[(a1)] if  $d=0$ (mod $3$)  , $D$ is very ample if and only
if $b+a\mu^-(E)
\geq 3$
	\item[(a2)]if  $d=1$ (mod $3$) , $D$ is very ample if $b+a\mu^-(E)
> 1$
	\item[(a3)]if  $d=2$ (mod $3$) , $D$ is very ample if
$b+a\mu^-(E) >\frac{4}{3}$
  \end{itemize}
\item[(b)] if $E$ is decomposable, then $D$ is very ample if and only if
$b+a\mu^-(E) \geq 3$
except when $E=E_1\oplus E_2$, with rk$E_1=1$, rk$E_2=2$, deg$E_2$ odd and
deg$E_1> \frac {degE_2}2.$ In the latter  case the condition is only
sufficient.
 \end{itemize}
 \end{introteo}
 Notice that the above theorem shows the existence, among others, of a
smooth threefold of
 degree 20 embedded in $\Bbb{P}^9$ as a fibration of Veronese surfaces
over an
 elliptic curve, choosing $a = 2, b= -1$ and $ d= 4.$

The authors would like to thank Enrique Arrondo and Antonio
Lanteri for their friendly advice.
The third author would like to thank the Department of Mathematics of
Oklahoma
State  University for the kind hospitality and warm support during the
final stages of
this work.
\section{General Results and Preliminaries}
\label{prelimsec}
\subsection{Notation}
\label{notation}
The notation used in this  work is mostly standard from Algebraic
Geometry. Good references are \cite{H} and \cite{gh}.
The ground field is always the field $\Bbb{ C}$ of complex
numbers. Unless otherwise stated all varieties are supposed to be
projective.
$\Bbb{P}^{n}$ denotes the n-dimensional complex projective space and
$\Bbb{C}^*$
the multiplicative group of non zero complex numbers.
 Given a projective n-dimensional variety $X$, ${\cal O}_X$
denotes its structure sheaf  and  $Pic(X)$  denotes the group of line
bundles over $X.$ Line bundles, vector bundles and Cartier
divisors are denoted by capital letters as $L, M, \dots.$ Locally free
sheaves of rank
one, line bundles and Cartier divisors are used interchangeably as
customary. Let $L, M \in Pic(X)$, let $E$ be a vector bundle of rank  $r$
on $X$, let
$\cal{F}$ be a coherent sheaf on $X$ and let $Y\subset X$ be a subvariety
of $X.$
Then the following notation is used:
\begin{enumerate}
\item[ ] $L C$ the intersection number of $L$ with a curve $C,$
\item[ ] $L^{n}$ the degree of $L,$
\item[ ] $|L|$ the complete linear system of effective divisors
associated with $L$,
\item[ ]$L_{|Y}$ the restriction of $L$ to  $Y,$
\item[ ] $L \sim M$ the linear equivalence of divisors
\item[ ] $L \equiv M$ the numerical equivalence of divisors
\item[ ] Num$(X)$ the group of line bundles on $X$ modulo the numerical
equivalence
\item[ ] $ E^*$ the dual of $E.$
\item[ ]  $\Bbb{P}(E)$ the projectivized bundle of $E$
\item[ ] $H^i(X, \cal{F})$ the $i^{th}$ cohomology vector space with
coefficient in ${\cal
F},$
\item[ ] $h^i(X,\cal{F})$ the dimension of $H^i(X, \cal{F}).$
 \end{enumerate}

If $C$ denotes a smooth projective curve of genus $ g$, and $E$ a vector
bundle
over $C$ of deg $E= c_1(E)= $d and rk $E=r$, we need the following
standard
definitions:
\begin{enumerate}
\item[ ] $E$ is $\it normalized$ if $h^0(E)\ne 0$ and $h^0(E \otimes L)=0$
for any invertible sheaf $L$ over $C$ with deg$L<0$.
\item[ ] $E$ has slope $\mu(E) = \frac{d}{r}$.
\item[ ] $E$ is $\it semistable$ if and only if for every proper subbundle
$S$,
$\mu(S) \leq \mu(E)$. It is $\it stable$ if and only if the equality is
strict.
\item[ ] The Harder-Narasimhan filtration of $E$ is the
unique filtration:
$$0=E_0\subset E_1\subset ....\subset E_s=E$$
 such that
$\frac{E_i}{E_{i-1}}$ is semistable for all $i$, and $\mu_i(E)=\mu
(\frac{E_i}{E_{i-1}})$ is
a strictly  decreasing function of $\it i$.
\end{enumerate}
We recall now some definitions from  \cite{bu} which we will use in the
following:
let  $0=E_0 \subset E_1 \subset ....\subset E_s=E$ be the Harder-Narasiman
filtration
of a vector bundle $E$ over $C$. Then
\begin{enumerate}
\item[]$\mu^-(E)=\mu_s(E)=\mu (\frac{E_s}{E_{s-1}})$
\item[]$\mu^+(E)=\mu_1(E)=\mu (E_1)$
\item[]or alternatively
\item[]$\mu^+(E)= $max $\{\mu(S) |0 \to S \to E \}$
\item[]$\mu^-(E)= $min $\{\mu(Q) |E \to Q \to 0 \}$.
\end{enumerate}
We have also $\mu^+(E) \geq \mu(E) \geq \mu^-(E)$ with equality if and
only if $E$ is
semistable.
In particular if $C$ is an elliptic curve, an indecomposable vector bundle
$E$ on $C$ is
semistable and hence  $\mu(E) = \mu^-(E)$.
Moreover if  $F,G$  are
indecomposable and hence semistable vector bundles on an elliptic curve
$C$ and $F
\to G$ is a non zero map, it follows that $\mu(F)\leq\mu(G)$.
\medskip
\subsection{General Results}
\medskip
Let $C$ be a smooth  projective curve of genus $g$, $E$ a vector
bundle of rank $r$, with $r \ge 2$, over $C$ and $\pi : X =\Bbb{P}(E)
 \to C $ the projective bundle
associated to $E$ with the natural projection $\pi$.
 With standard  notations denote
with $\cal {T} = \cal {O}_{\Bbb{P}(E)} (1) $ the tautological sheaf and
with  $\cal {F}_P=
\pi^*\cal {O}_{C}(P) $  the line bundle associated with the fiber over
$P\in C.$ Let
$T$ and
$f$ denote the numerical classes respectively of $\cal T$ and $\cal
{F}_P$.

  Let $D\sim a\cal{T} + \pi^*B$, with  $ a\in \Bbb{Z}$, $B\in Pic(C)$
and deg$B = b$, then $  D \equiv aT+bf .$  Moreover  $\pi_*D = S^{a}(E)
\otimes
\cal {O}_{C}(B) $ and hence $\mu^-(\pi_*D)=a\mu^-(E) +b$ (see \cite{bu}).

 Regarding  the ampleness, the global generation, and
the normal generation  of $D$, the following criteria are known:
 \begin{Thm}[Miyaoka \cite{Miyao3}]
\label{miyaoteo}
 Let $E$ be a vector bundle over a smooth
projective curve $C$ of genus $g$, and $X =\Bbb{P}(E)$ . If $D\equiv
aT+bf$ is a line
bundle over $X$, then $D$ is ample if and only if $a>0$ and $b+a \mu^-(E)
>0$.
\end {Thm}

\begin{Prop}[Gushel  \cite{gu2},  proposition 3.3]
\label{guongg}
Let $D \sim a\cal{T}+\pi^*B$ where $a > 0$ and $B \in$ Pic($C$), be a
divisor on a
projective bundle $\pi : X =\Bbb{P}(E) \to C $ . Then:
\begin{itemize}
\item[i)]if $a=1$, the bundle $\pi_* (D)$  is
generated by global sections  if and only if the divisor $D$ is
\item[ii)]if $a \geq 2$, and the vector bundle $\pi_*(D)$
 is
generated by global sections, then also the divisor $D$ is.
\end {itemize}
\end {Prop}

\begin {Lem}[Gushel \cite{gu}, Proposition 3.2]
\label{ggforindec}
 Let $E$ be an indecomposable vector bundle  over an elliptic curve $C$.
$E$  is
globally generated if and only if deg$E > $rank$E$.
\end{Lem}
\begin{Lem}
\label{buongg}
(see e.g. \cite{bu}, lemma 1.12)
Let $E$ be a vector bundle over $C$ of genus $g$.
\begin{itemize}
\item[i)]if $\mu^-(E) > 2g-2$ then $h^1(C,E)=0$
\item[ii)]if $\mu^-(E) > 2g-1$ then $E$ is generated by global sections.
\end{itemize}
\end{Lem}
For the following theorem we need a definition:
\begin{Def}[Butler, \cite{bu}]
Let $E$ be a vector bundle over a variety $Y$, and let $\pi: X =\Bbb{P}(E)
\to Y$ be
the natural projection. A coherent sheaf $\cal F$  over $X$ is said to be
$t \pi-regular$ if, for all  $i>0$, $$ \cal {R}^{i}\pi_{*}(\cal {F}(t-i))
=0. $$
\end{Def}
\begin{Thm}[Butler,\cite{bu}]
 Let $E$ be a vector bundle over a smooth
projective curve $C$ of genus $g$,  and $X =\Bbb{P}(E)$ . If $D$ is a
$(-1){\pi}-$ regular line bundle over X, with $\mu^-(\pi_*D) >2g $, then $
D$ is
normally generated.
\label{butlerteo}
\end{Thm}

\begin{rem}
\label{criteriodelbutler}
 Let $D$ be a divisor of $X =\Bbb{P}(E)$, with $E$ vector bundle on a
smooth
projective curve of genus $g.$ As  $h^{i}(\cal {F}_p,D_
{|\cal{F}_p}(-1-i))=0$ for $i\ge1$
, the $(-1)\pi $- regularity of $D$ is satisfied, hence the condition
$a\mu^-(E) +b
>2g$ implies that $D$ is normally generated.
If a line bundle $D$ on a projective variety $X$ is ample and
normally generated it is very ample. Hence from Theorem \ref{miyaoteo} and
\ref{butlerteo}  we get that $D$ is very ample on $ X =\Bbb{P}(E) $ if
\begin{equation}
\label{condizionedelbutler}
b+a\mu^-(E) > 2g.
\end{equation}  Hence, if
$g=1$, the very ampleness of $D \equiv aT+bf$  is an
open problem only in the range
\begin{equation}
\label{range}
0<b+a\mu^-(E) \le2.
\end{equation}
\end{rem}
\medskip
 \subsection{Preliminaries}
\medskip
The following result is standard from the theory of vector bundles
(see \cite{H}):
\begin {Lem}
\label{degofnormed}
Let $E$ be an indecomposable vector bundle of rank $r$ on an elliptic
curve.
If $E$ is normalized then $0 \leq deg E \leq r-1$.
\end {Lem}

\begin{Lem}
Let $E=\bigoplus_{i=1}^{n} E_i$ be a decomposable vector bundle over an
elliptic
curve $C$, with $E_i$ indecomposable vector bundles. Then $\mu^-(E)=$ min$
 \,\mu(E_i)$.
\label{mimenodec}
\end{Lem}
\begin{pf}
For the proof we need the following three claims.
\begin{claim}
 Let $E=\bigoplus_{i} E_i$ be as above,
then $\mu(E) \ge$ min$\,\mu(E_i)$.
\label{claim1}
\end{claim}
\begin{pf} Let us denote by $r=rk(E)$ $r_i=rk(E_i)$ $d=deg(E)$,
$d_i=deg(E_i)$.
 Let us
consider the vectors $\underline{v}_i$ in $\Bbb{R}^2$ whose coordinates
are
$(r_i,d_i)$ and the vector $\underline{v} = \sum_{i} \underline{v}_i$. Let
$\alpha_i$
be the angle between the $r-$axis and $\underline {v}_i$. Let $\alpha$ be
the angle
between the $r-$axis and $\underline {v}$. It is   $\mu(E) = \frac{d}{r}=$
tg$(\alpha) \geq $ min$_i$ tg$(\alpha_i) = $ min$_i$ $(\frac{d_i}{r_i}) =
$ min$_i$ $
\mu(E_i)$.
\end{pf}
 \begin{claim} Let  $E=\bigoplus_{i} E_i$ be as above, and
$\mu(E_i)=\frac{d_i}{r_i}=h \in \Bbb {Q},$  for all  $i$. Then  $\,
\mu^-(E) = h$.
\label{claim2}
\end{claim}
\begin{pf}
Notice that  under this hypothesis  $\mu(E)=h$. Moreover, by definition,
it is
$\mu^-(E)
=$ min $\{\mu(Q)\, | E \to Q \to 0 \}$. If $Q$ is decomposable in the
direct sum of
indecomposable vector bundles $ Q_k$, the existence of a surjective map
$ E \to Q \to0 $ implies the existence of surjective maps $ E \to Q_k \to
0 $  for all  $k$
and consequently from Claim \ref{claim1},
 $\mu^-(E) =$ min $\{\mu(Q)\, | E \to Q \to 0,$ and $ Q$ indecomposable
\}.

Now let $Q_o$ be an
indecomposable vector bundle which realizes the minimum, i.e.
$\mu(Q_o)= \mu^-(E).$ From $\oplus_i  E_i \to Q_o  \to 0 $  it follows
that there exists
at least an index
${i_0}$  such that the map $E_{i_0} \to Q$ is not zero and
$\mu(E_{i_0})\leq\mu(Q_o).$
Therefore it is  $h\leq\mu^-(E)$. As $h\geq\mu^-(E)$, the Claim is proved.
\end{pf}
\begin{claim} Let  $E=\bigoplus_{i}E_i$ be as in Claim (2). Then E  is
semistable.
\label{claim3}
\end{claim}
\begin{pf}
It is enough to prove that for any $ S$ vector bundle on $C$ such that
there exists
a map $0 \to S \to E$ then $\mu(S) \leq \mu(E) =h$. If we consider the
dual map
 $E^* \to S^* \to 0$ we have $\mu(S^*) \geq \mu^-(E^*)
=\mu^-(\bigoplus_{i}
 E{_i}{^*})$ and, as $\mu(E^*{_i}) =- \frac{d_i}{r_i} = -h$, from Claim
\ref{claim2} applied
to
$E^*$ we have $\mu^-(E^*)= -h$. Hence $\mu(S^*)=-\mu(S) \geq -h$ and
$\mu(S) \leq h$.
\end{pf}

The Lemma can now be proved.

Let  $E=\bigoplus_{i} E_i$ be as in the hypothesis of Lemma, and denote by
$\mu_i =\mu(E_i)$. We can choose an ordering  such that $E=E_1 \oplus E_2
\oplus E_3 .....$ and $\mu_1\geq \mu_2 \geq \mu_3 ...$.
Let  $E=\bigoplus_{k=1}^{s} A_k$ be a new decomposition of $E$ such that
each $A_k$
is an indecomposable vector bundle or a sum of indecomposable vector
bundles
$E_i$ with the same $\mu_i$. In this way we get a strictly decreasing
sequence
$\mu(A_1)>\mu(A_2)>...>\mu(A_s)$, and by claim (3) each $A_k$ is
semistable.
Moreover the sequence $ 0\subset F_1 \subset F_2 \subset ... \subset F_s
=E$
with $F_i = A_1 \oplus A_2 \oplus ...\oplus A_i$ with $ 1\leq i  \leq s$,
is the
Harder-Narasimhan filtration of $E$ because the sequence of the slopes
$\mu(\frac
{F_i}{F_{i-1}}) = \mu(A_i)$ is  strictly decreasing and each $\frac
{F_i}{F_{i-1}}=A_i$ is
semistable  for all  $i=1...s$.
Hence we get $ \mu^-(E) = \mu(\frac {E_s}{E_{s-1}}) = \mu(A_s) =$ min $
\mu(E_i).$
\end {pf}

\begin{Lem}
Let $D\sim a\cal{T} + \pi^*B$ be a line bundle
in $X =\Bbb{P}(E)$ over a curve $C$ of genus $g=1$,
with $B\in Pic(C)$ , $a \geq 1$ and deg$B = b $ .
\begin{itemize}
\item[i)]If $a=1$ $D$ is globally generated  if and only if  $b+ \mu^-(E)
> 1$
\item[ii)]If $a\geq2$ $D$ is globally generated  if  $b+a\mu^-(E) > 1$
\end{itemize}
\end{Lem}
\begin{pf}
To prove   ii) it is sufficient to apply  Proposition \ref{guongg} and
Lemma
\ref{buongg}. To prove  i)  notice that if $E$ is indecomposable it is
enough to apply
Proposition \ref{guongg} and Lemma \ref{ggforindec} , observing that an
indecomposable vector bundle $E$ over an elliptic curve is semistable and
hence
$\mu^-(E)=\mu(E).$
Let now  $E$ be decomposable and hence $E \otimes B$
decomposable over $C$. In particular let $ E \otimes B=
\bigoplus_{q=1}^{s}
A_q $ be a decomposition of $ E \otimes B$ in indecomposable vector
bundles
$A_q$ over $C$. By Lemma \ref{ggforindec}
 every $A_q$, for $q=1...s$ is
globally generated if and only if deg $A_q >$rk $A_q$, i.e. if and only if
$\mu(A_q)>1,$
for all  $q$. From Lemma \ref{ggforindec}, Lemma\ref{mimenodec} and
Proposition
\ref{guongg}  we get the following chain of equivalences which conclude
the proof:
$\mu^- (E)+b > 1\Leftrightarrow \mu^-( E \otimes B)=$ min$_q  \mu(A_q) >1
\Leftrightarrow \mu(A_q) >1$  for all  $q$  $\Leftrightarrow A_q$ is
globally
generated
 for all  $q$ $ \Leftrightarrow \pi_*D$ is globally generated on $C
\Leftrightarrow  D $
is globally generated on $X$
\end{pf}
The above Lemma is partially contained in \cite[Prop.
3.3]{gu}. Unfortunately  the proof presented there is based on
\cite[Prop.1.1 (iv)]{gu},
which is not correct, as the following counterexample shows.
Let $E$ be an indecomposable  vector bundle over an elliptic curve with
deg $E
= 1$ and rank $E = 2.$ Then $2 \cal{T} = \cal{O}_{\Bbb{P}(E)}(2)$ is
generated by global sections, according to \cite[Prop. 8.5.8]{BESO}. On
the
other hand let $\pi_*(2\cal{T}) = S^2 E = \bigoplus_q A_q,$ where $A_q$ is
indecomposable for all $q.$ Then $S^2E$ is
generated by global sections if and only if $A_q$ is such, for all $q.$
>From Lemma \ref{ggforindec} it follows that $S^2E$ is globally generated
if and
only if $\mu(A_q) > 1$ for all $q,$ i.e. if and only if $\mu^-(S^2E) > 1,$
i.e. if and only if $ 2 \mu^-(E) = 2 \mu(E) > 1$ which is false.

\medskip
If we consider an indecomposable vector bundle of degree  $d=0$, we have
the
following proposition. It is contained in
\cite[Theorem 3.9]{gu2}, but we prefer to give here a simpler proof.
\begin{Prop}
\label{vadeg0modr}
Let $E$ be an indecomposable rank r vector bundle over an elliptic curve
$C$ with
deg$E=0$ (mod $r$), and let $D\equiv aT+bf$ be a line bundle on $ X
=\Bbb{P}(E) $.
Then $D$ is very ample on $X$ if and only if $b+a\mu^-(E) = b+a\mu(E)  \ge
3$.
\end {Prop}
\begin{pf}
It is enough to consider  the case in which $E$ is normalized , as
if $E$ is not normalized we can consider its normalization
$\bar {E}=E \otimes L$ with deg$L=l$. If $D \equiv aT+bf$  in
Num$\Bbb{P}(E)$ then
  in Num$\Bbb{P}(\bar{E})$ we get
\begin{equation}
\label{contonorm}
D \equiv a{\bar{T}}+(b-al)\bar{f}\thinspace,\ \ \ \ \bar{d}=
deg\bar{E}=d+rl,\
\ \ \  \mu(\bar{E}) = \mu(E) + al.
\end{equation}
 Let $E$ be  normalized, hence $d=0$ and $E =
F_r$ in  the notation of \cite {At} (recall that $F_1 = \cal{O}_C$).
According to  (\ref{range})
the only
cases to be
considered are
$b=1$ and
$b=2$ and hence
$D\equiv aT+f$  or
$D\equiv aT+2f$.  We want to show that in both these cases $D$ is not very
ample.
Assume the contrary and  proceed by induction on $r$. Let  $r=2.$ As $D
T =
1$ or
$2$, the smooth elliptic curve $\Gamma$, which is the only element of
$|\cal T|$, is embedded by  $\phi_{|D|}$ as a line or a conic which is a
contradiction.
Assume now  the proposition true for $F_{r-1}$  and recall that  there is
a
short exact sequence (see \cite{At} pag 432)
\begin{equation}
\label{succdegliFr}
0 \to \cal{O}_C \to F_r
\to F_{r-1}
\to 0.\end{equation}
Let  $T'= T_{|Y}$ and $f' = f_{|Y}$ the generators of Num$(Y)$ where
$ Y =\Bbb{P}(F_{r-1}) \subset X =\Bbb{P}(E) $. If $D$ is very ample,
$D_{|Y}$ is very
ample too; but $D_{|Y} \equiv aT'+bf'$ and it is not very ample by
induction
hypothesis. Hence $D$ is very ample if and only if $b \geq 3$.
\end {pf}
The following Lemma, which gives a sufficient condition for the very
ampleness
of a divisor $D$ on $X =\Bbb{P}(E),$ will be needed later on.
\begin {Lem}
\label{lemmadiEnrique}
Let $E$ be a rank $r$ vector bundle over a curve $C$ and let $D\equiv
aT+bf$ be a line bundle on
$ X =\Bbb{P}(E) $,with $a \geq 1$. If $ \pi_*D$ is a very ample vector
bundle on the
curve $C$, then $D$ is very ample on $\Bbb{P}(E)$. Moreover if $a=1$, $D$
is very ample on $X$
if and only if $\pi_*D$ is very ample on $C.$
\end {Lem}
\begin {pf} We give only a sketch of the proof.  A divisor  $D\equiv
aT+bf$ on $X$ defines a map
$\varphi_{|D|}$ in a suitable projective space such that $X'=
\varphi_{|D|}(X)$ is a bundle on $C$
 whose fibers are the
Veronese embedding of the fibers of $X =\Bbb{P}(E)$. Moreover each fiber
of $X'$ is embedded in
a fiber of the projective bundle $\Bbb{P}(S^{a}(E) \otimes \cal
{O}_{C}(B))$.It follows that the
very ampleness
 of $S^{a}(E) \otimes \cal {O}_{C}(B)$ and hence of its tautological
bundle implies that the map
$\varphi_{|D|}$ gives an embedding and hence that $D$ is very ample.
The case $a=1$ follows immediatly from the above considerations.
\end {pf}
\medskip
\subsection{The case $a=1$}
\medskip
We want to investigate the very ampleness of $D \equiv aT+bf$ in
dependence of $a$
and $b$. As we have remarked at the end of section 2.2, the problem is
open only
when $0 < b+a\mu^-(E)  \leq 2$. Let us begin with the case $a=1$.
 In this case we have the following theorem:
\begin{Thm}[Gushel,\cite{gu} theorem 4.3]
 Let $D \sim \cal{T}+ \pi^*B$ be a divisor on
$\Bbb{P}(E) $, where $E$ is an indecomposable and normalized vector bundle
of rank $r$
over an elliptic curve $C$. If $b= $deg$B$, the divisor $D$ is very ample
if and only if:
\begin{enumerate}
\item[i)] $b \geq 3$ if  deg$E=0$
\item[ii)] $b \geq 2$ if $0< $deg$E < r$.
\end{enumerate}
\end{Thm}
Now it is easy to prove the following (see (\ref{contonorm})):
\begin{Prop}
\label{vaindeca=1}
In the above assumptions and notations, if $E$ is indecomposable
but not normalized, it follows that $D$ is very ample if and only if the
following
conditions hold :
\begin{enumerate}
\item[i)] $b + \mu(E) \geq 3$ if $d=0$ (mod $r$)  .
\item[ii)] $b + \mu(E) \geq 2$ otherwise.
\end{enumerate}
\end{Prop}
\medskip
\begin{rem} The previous results consider the case in which $E$ is
indecomposable.  If $E$ is decomposable, by Lemma 2.12, we can argue as
follows: firstly
in this case, as $a=1$,
$D$ is very ample if and only if $\pi_*(D)$ is very ample.  Secondly we
have
$D \sim \cal
T + \pi^{*}B$, $\pi_*(D) \simeq E \otimes \cal {O}_{C}(B) = \bigoplus E_j
\otimes
\cal {O}_{C}(B)$, with $E_j$ indecomposable vector bundles .  Moreover $E
\otimes \cal
{O}_{C}(B) $ is very ample if and only if every $E_j \otimes \cal
{O}_{C}(B) $ is very
ample.  Let deg$E_j = d_j$ and rk$E_j = r_j$, and assume that $d_j= 0$
(mod
$r_j$), possibly only for $j = 1...t.$ Then  $D$ is very ample if and only
if
$ b+ \frac{d_j}{r_j} \geq 3$ for $j=1...t$, and  $ b+ \frac{d_j}{r_j} \geq
2$  for the
remaining
$j$'s, by Proposition \ref{vaindeca=1}.
\end{rem}

 \medskip
Having dealt above with the case $a=1$, from now on the blanket assumption
$a\ge
2$ will be in effect.

\section{Rank 2}
Let $E$ be a rank $2$ vector bundle on an elliptic curve $C$ and let
$D\equiv aT+bf$
be a line bundle on  $ X =\Bbb{P}(E)$. Assume that $E$ is indecomposable.
If $E$
is normalized then deg$E = 0,1$ by Lemma \ref{degofnormed}. If deg$E=0,$
from
Proposition \ref{vadeg0modr} it follows  that $D$ is very ample if and
only if $b \geq
3$. If deg$E=1$, necessary and
sufficient conditions for the very ampleness of $D$ are given by the
following
Theorem, reformulated under our assumption that $a\ge 2.$
\begin{Thm}[Biancofiore - Livorni, \cite{bi-li3},Theo 6.3]
Let $D \sim a\cal{T}+ \pi^*B$ be a divisor on $\Bbb{P}(E) $, where $E$ is
an
indecomposable normalized
vector bundle of rank $2$ and degree $1$ over an elliptic curve $C$. If
$b= $deg$B$,
the divisor $D$ is very ample if and only if $b + \frac{a}{2} > 1$.
 \end{Thm}
The following Proposition can now be easily proved (see
(\ref{contonorm})).
\begin{Prop}
\label{rangodue}
In the above hypothesis, if $E$ is indecomposable but not
normalized,  $D$ is very ample if and only
if the following conditions hold :
\begin{enumerate}
\item[ii)] $b + a\mu^{-}(E) \geq 3$ if $d=0$ (mod $2$)  .
\item[i)] $b + a\mu^{-}(E) > 1$  if $d=1$ (mod $2$).
\end{enumerate}
\end{Prop}
The case $E$ decomposable is treated by  the following Theorem.
\begin{Thm}
Let $D \sim a\cal{T}+ \pi^*B$ be a divisor on $\Bbb{P}(E) $, where $E$ is
a decomposable
vector bundle of rank $2$ over an elliptic curve $C$, $b =$ deg $B$. The
divisor $D$ is very ample if and only if $b + a\mu^{-}(E) \geq 3.$
 \end{Thm}
\begin{pf}
To prove   the sufficient condition let $E$ be decomposable as
$H\bigoplus G$ where  $H$ and $G$ are line bundles on $C$ with deg$H=h
\geq$
deg $G=g$. By Lemma \ref{mimenodec}  it is  $\mu^{-}(E) = g$. By Lemma
\ref{lemmadiEnrique}, a sufficient condition for the very ampleness of $D$
on $X$ is
that
$\pi_*(D) $ is very ample as  a vector bundle on $C$. In our hypothesis

$$\pi_*(D) =  S^a(E) \otimes \cal {O}_C(B) =
\bigoplus_{q=0}^{a}H^{\otimes q} \otimes G^{\otimes a-q} \otimes B.$$
 Now $\pi_{*}(D) $ is very ample  if each element of its decomposition has
degree
$\geq 3$, i.e. if $qh + (a-q)g +b \geq 3$,  for all  $q=0,...a$. As the
minimum of $qh +
(a-q)g +b$ is realized for $q=0$, $\pi_{*}(D) $ is very ample if and only
if $ag +b = b+a \mu^{-}(E)
\geq 3$. This condition is also necessary for the very ampleness of $D$.
Indeed
 the projective bundle $\Bbb{P}(G) $, by the exact sequence
 $$0 \to H \to E \to G \to 0$$ \cite{gu}, Proposition 1.1,  gives an
elliptic
curve $\Gamma$ on
$X,$  $\Gamma \in | \cal{T}+ \pi^{*}(H^{*})| .$ Notice that
$h^0(X, \cal{T}+ \pi^{*}(H^{*}))>0.$ If $D$ is very ample it must be $D
\Gamma =
(aT+bf)(T-hf)= ag+b= b+a
\mu^{-}(E) \geq 3$.
\end{pf}
\medskip

\section{Rank 3}
\label{rango3sec}
In this section and in the next one, we will prove the very ampleness of a
divisor on a
smooth variety following a classical  method,  based on the following
lemmata.
\begin{Lem}
\label{reductiontoS}
Let $X$ be a smooth variety, $D$ a divisor in $Pic(X)$ and let $A$ be
another element
of
$Pic(X)$, such that $h^1(X, \cal{O}_X(D-A))=0$. If,  for each pair of
points $R,Q \in X$
(possibly  infinitely near) it is possible to find a smooth element $S \in
|A|$ containing
$R$ and $Q$, then $D$ is very ample on $X$ if and only if $D_{|S}$ is very
ample on $S$.

\end{Lem}
\begin{pf}
If $D$ is very ample then obviously $D_{|S}$ is very ample. On the other
hand, pick
any two points $R,Q \in X$ (distinct or infinitely near) and choose $S \in
|A|$ such
that  $R,Q \in S$. As $D_{|S}$ is very ample there exist sections of
$D_{|S}$ separating
$R$ and $Q$. Now look at the following exact sequence
 $$0 \to \cal{O}_X(D-A) \to \cal{O}_X(D) \to \cal{O}_S(D_{|S}) \to 0.$$
 From the
assumptions above
 we get that the map $H^0(X,\cal{O}_X(D)) \to  H^0(S,\cal{O}_S(D_{|S}))$
is surjective and hence $D$ is very ample on $X$ if and only if $D_{|S}$
is
very ample on $S$.
\end{pf}
\begin{Lem}
\label{T+F}
Let $E$ be an ample vector bundle over an elliptic curve $C$ such that
deg~$E<
$~rk~$E$.
 Let $ X =\Bbb{P}(E) $, let
$P$ be a fixed point of
$C$,
$D\sim a\cal{T} +\pi^*B$ and
$A = \cal{T}+\cal{F_P}$ be   line bundles on  $X$, with $b =$ deg $B$ ,
$b-1+(a-1)\mu^-(E) >0$ and
 $h^0(X,\cal{O}_X(A)) \geq $ deg $E +3$. Then the hypothesis of Lemma
\ref{reductiontoS} are satisfied for $A$.
 \end{Lem}
\begin{pf}
If $ (a-1)\mu^-(E)+b-1>0$ then by Lemma \ref{buongg} it is   $h^1(X,
\cal{O}_X(D-A))=0.$ Moreover being  $h^0(X, \cal{O}_X(A)) \geq  $ deg $E
+3,$
for each pair of points $R, Q \in X$ there exists a linear subsystem $\cal
L \subset |A|$,
with  dim$\cal L \ge$  deg $E$ , all the elements of which contain $R$ and
$Q$.
Moreover in
$\cal L$ there is at least one smooth element $S.$  In fact, assume that
all the
elements of  $\cal L$ are singular. Note that any singular element of
$\cal L$ must
be reducible as $\Gamma \cup \cal{F}_P$, with $\Gamma \in |\cal{T}|$,
$\Gamma$
smooth, because we have that any divisor numerically equivalent to $T-f$
is not effective as
deg $E< $ rk $E$. As $h^0(X, \cal{O}_X(\cal{F}_p)) =1,$  for all  $P \in
C$, by Bertini's theorem  all the elements of $\cal L$ are singular only
if $\cal{F_P}$ is
fixed and  $\Gamma$ varies in a subsystem of $|\cal{T}|$ of dimension deg
$E$.  This is
impossible  as   $h^0(X,\cal{O}_X(\cal{T}))= $ deg $E$.
\end{pf}
\begin{Lem}
Let $E$ be an ample vector bundle over an elliptic curve $C$ such that
deg$E<$rk$E$.
 Let $ X =\Bbb{P}(E) $, and let $D \sim a\cal{T}+\pi^*B$  be  a
line bundle on  $X$, with $b= $ deg $B$ , $b+(a-1)\mu^-(E) >0$ and
 $h^0(X,\cal{O}_X(\cal{T})) \geq 3.$ Then the hypothesis of Lemma
\ref{reductiontoS}
are satisfied, with $A = \cal T$.
\label{T}
 \end{Lem}
\begin{pf}
If $ (a-1)\mu^-(E)+b>0$ then Lemma \ref{buongg} gives $h^1(X,
\cal{O}_X(D-A))=0.$ Moreover as each element of $|\cal{T}|$ is smooth,
because we have that any
divisor numerically equivalent to $T-f$ is not effective as
deg$E<$rk$E$, the condition
$h^0(X,\cal{O}_X(\cal{T})) \geq 3$ shows that it is possible to find  a
smooth element $S
\in |\cal{T}|$ containing  each fixed pair of points $R,Q \in X$.
\end{pf}

The following Lemma will be very useful to obtain the vanishing condition
required
by Lemma \ref{reductiontoS} in many borderline cases. The notation used
here is
the classical notation used by Atiyah in \cite{At}.
\begin{Lem}
\label{AA}
Let $E$ be an indecomposable vector bundle over an elliptic curve $C$ with
rank $E =
r$ and deg $E = d.$ Let
$X=\Bbb{P}(E)$ and let $\pi : X \to C$ be the natural projection.
Let $D =a\cal{T}+ \pi^*(B)$ for a line bundle $B$ with
deg $B = b$ and let $A=\cal{T}
+\pi^*(\cal{O}_C(P))$  where $P$ is a point in $C.$

If  $\frac{(a-1) d}{r} + b -1=0,$ it is possible to choose $P\in C$ such
that $h^1(X, D - A) = 0.$
\end{Lem}
\begin{pf}
It is enough to show that  $h^1(C, S^{a-1} E \otimes B \otimes
\cal{O}_C(-P)) = 0.$

Since deg $S^{a-1} E \otimes B \otimes \cal{O}_C(-P) = 0$ by Riemann Roch
it is
enough to show that  $h^0(S^{a-1} E \otimes B \otimes \cal{O}_C(-P)) = 0.$
Because $S^{a-1}(E)$ is a direct summand of $E^{\otimes (a-1)}$ it is
enough to show
that $h^0( E^{\otimes (a-1)} \otimes B~\otimes \cal{O}_C(-P)) = 0.$

Let $h =$ gcd $(d, r).$ Then by \cite{At} Lemma 24 and 26 it is
\begin{equation}
\label{EconFh}
E = E' \otimes F_h
\end{equation}
where $d'=$ deg $E' = \frac{d}{h}$ and $r' = $ rank $E' = \frac{r}{h}$ so
that gcd$(d',r')
= 1,$ $F_h$ is as in \cite{At} Theorem 5 and  $E'$ is indecomposable  .
Being  $r'$ and $d'$ relatively prime, the condition $\frac{(a-1) d'}{r'}
+ b
-1=0$ shows that $r'$ divides $(a-1).$ Therefore following \cite{H2}
Proposition 1.4
it follows that
\begin{equation}
\label{EprimoconFri}
 E^{' \otimes (a-1)} = \bigoplus_i(F_{r_i} \otimes L_i)
\end{equation}
Therefore putting (\ref{EconFh}) and (\ref{EprimoconFri}) together we get
$$ E^{\otimes(a-1)} =( E' \otimes F_h)^{\otimes (a-1)} =
\bigoplus_i(F_{r_i} \otimes
L_i) \otimes F_h^{\otimes (a-1)}.$$

Theorem 8 in \cite{At} shows that tensor powers of $F_l$ 's are direct
sums of $F_k$'s
so we conclude that
$$ E^{(a-1)} = \bigoplus_j (F_{r_j} \otimes L_j).$$
It is then enough to show that for all $j$ it is $h^0(F_{r_j} \otimes L_j
\otimes B
\otimes \cal{O}_C(-P)) = 0.$

Let $\cal{L}_{j,P} = L_j \otimes B\otimes \cal{O}_C(-P) .$
Recall that the $F_{r_j}$ are obtained as successive extensions of each
other by
$\cal{O}_C,$
i.e. for every $r$  we have the sequence (\ref{succdegliFr}) (see proof of
Prop .2.11).

This shows that it is $h^0(F_{r_j} \otimes \cal{L}_{j,P})  = 0$  unless
$\cal{L}_{j,P} =
\cal{O}_C.$
It is then enough to choose a point $P$ such that $L_j \otimes B \otimes
\cal{O}_C(-P)
\neq \cal{O}_C$ for all $j.$ Since $B$ is  a fixed line bundle and $j$
runs over a finite
set, a $P$ that works for all $j$ can certainly be found.
\end{pf}

 The following Theorem collects our results for  the case rk $E = 3.$
\begin{Thm}
\label{Thmrango3}
Let $E$ be a rank $3$ vector bundle on an elliptic curve $C$ and let
$D\equiv aT+bf$
be a line bundle on  $ X =\Bbb{P}(E).$
\begin{itemize}
\item[(a)] If $E$ is indecomposable then
	\begin{itemize}
	\item[(a1)] if  $d=0$ (mod $3$)  , $D$ is very ample if and only
if $b+a\mu^-(E)
\geq 3$
	\item[(a2)]if  $d=1$ (mod $3$) , $D$ is very ample if $b+a\mu^-(E)
> 1$
	\item[(a3)]if  $d=2$ (mod $3$) , $D$ is very ample if
$b+a\mu^-(E) >\frac{4}{3}$
  \end{itemize}
\item[(b)] if $E$ is decomposable, then $D$ is very ample if and only if
$b+a\mu^-(E) \geq 3$
except when $E=E_1\oplus E_2$, with rk$E_1=1$, rk$E_2=2$, deg$E_2$ odd and
deg$E_1> \frac {degE_2}2.$ In the latter  case the condition is only
sufficient.
 \end{itemize}
 \end{Thm}

\begin{pf}
 Firstly we consider the case $E$ indecomposable and normalized. By Lemma
\ref{degofnormed} and Proposition \ref{vadeg0modr}  only the cases
$d=1$ and
$d=2$ need to be considered.
\begin{case}
$d=1$.
\end{case}
 Let $A$
be as in Lemma \ref{T+F} and  notice that $h^0(A) =$ deg $E + 3.$  By
Remark
\ref{criteriodelbutler}.
 we can assume $b+\frac{a}{3}= 1+\frac {\varepsilon}{3},$ where
$\varepsilon = 1,2,3.$   If
$\varepsilon =2,3$ then  $ (a-1)\mu^-(E)+b-1 = b+
\frac {a}{3} -
\frac {4}{3} >0.$
 If $\varepsilon = 1$
Lemma \ref{AA} still allows the use of Lemma \ref{reductiontoS}. Therefore
Lemma
\ref{reductiontoS} and \ref{T+F} can be applied. Let
$S$ be  a smooth element in $|\cal{T}+\cal{F_P}|.$ It is enough to check
when $D_{|S}$ is
very ample. By \cite{gu}, prop 1.1,
$S=\Bbb{P}(E')$ where $E'$ is a rank $2$ vector bundle on the curve $C$,
with
deg$E'=2$, defined by the sequence
\begin{equation} 0 \to \cal {O_C}(-P) \to E \to E' \to 0.
\label{secondasucc}
\end{equation}

As Num($S$) is generated by $T'$ and $f'$, with $T'=T_{|S}$ and
$f'=f_{|S}$, $D_{|S}
\equiv aT'+bf'$ and moreover by section $3$, $D_{|S}$ is very ample in our
range
if and only if        $b+ a\mu^-(E') \geq 3$. If $E'$ is indecomposable
then
$\mu^-(E')=\mu(E')=1$. If $E'$ is decomposable, i.e. $E'= H \oplus G $,
both $\mu(H)
\geq \mu(E) = \frac{1}{3}$  and   $\mu(G) \geq \mu(E) = \frac{1}{3}$
because from (\ref{secondasucc})
there exist  non zero surjective morphisms  $E \to H$ and   $E \to G.$

Hence deg$H$= deg$G$ =$1$, $\mu(H)=\mu(G)= 1$ and $\mu^-(E')=1.$ In every
case
the condition  $b+ a\mu^-(E') \geq 3$ is satisfied in the  range under
consideration.
\begin{case}
 $d=2$.
\end{case}
Let $A$ and $S$ be as in Lemma \ref{reductiontoS} and
\ref{T+F}. In this case it is
$h^0(X,A) = 5 \geq $ deg $E + 3$. By Remark \ref{criteriodelbutler} we can
assume
$b+a\frac{2}{3}=1+\frac{\varepsilon}{3}$ with $\varepsilon=1,2,3$. If
$\varepsilon=3$ the hypothesis of Lemma \ref{T+F} are satisfied. If
$\varepsilon=2$  Lemma \ref{AA} still allows the use of Lemma
\ref{reductiontoS}. Therefore it  suffices to investigate the very
ampleness of
$D_{|S}.$
Let $S=\Bbb{P}(E')$ with $E'$ a rank $2$ vector bundle on $C$, with
deg$E'=3$ defined again by (\ref{secondasucc}). If $E'$ is indecomposable,
$D_{|S}
\equiv aT'+bf'$ is very ample if and only if $b+a\mu^-(E') > 1$, i.e.
$b+\frac{3}{2} a >1$
i.e. for all $D$ in the range under consideration.

If $E'$ is decomposable then  $E'= H \oplus G $, with deg$H$ and deg$G$
$\geq
\mu(E)=\frac{2}{3}.$ Hence  we can assume  deg$H=1$ and deg$G=2.$ By Lemma
\ref{mimenodec}, $ \mu^-(E)=1$ and the very ampleness condition is $b+a
\geq 3$
which is satisfied by every  $D$  in the  range under consideration.

Notice  that in the case $b+\frac{2}{3} a = \frac{4}{3}$, i.e.
$\varepsilon=1$,
 a very ampleness result for all $D$ in our  range cannot be expected. For
example,
$D\sim 2\cal{T}$ is not very ample as $D_{|Y}$ is not very ample for each
smooth
surface
$Y\in |\cal{T}|$ by section $3$.

If  $E$ is indecomposable but not normalized, the result is obtained   by
 similar computations (see (\ref{contonorm})).

To prove (b), let  now $E$ be  decomposable.
Firstly we prove the sufficient condition. By
Proposition \ref{guongg} it suffices to prove that $\pi_{*}(D)$ is very
ample.
Let us consider $\pi_{*}(D) = S^{a}(E) \otimes \cal {O}_C(B) =
\bigoplus_{q}A_q$ where
$A_q$ is an  indecomposable vector bundle on $C$ for all $q.$
 $ S^{a}(E) \otimes \cal {O}_C(B)$ is very ample if and only if $A_q$ is
very ample
 for all $ q$ i.e. if $\mu(A_q)\geq 3 $   for all  $q.$
This condition is satisfied if   min$_q \mu(A_q) = \mu^-( S^{a}(E) \otimes
\cal
{O}_C(B))=b+a\mu^-(E)
\geq 3$ which is what we wanted to show.

To prove the necessary condition,   two cases will be considered:

i) $E$ is sum of three line bundles,  $E = W \oplus G \oplus H$
respectively of degrees
$w\leq g \leq h$. By lemma \ref{mimenodec}  $\mu^-(E) = w$,  and $d=
$deg$E=w+g+h.$ From \cite{gu}, Prop 1.1 it follows that
$\Bbb{P}(G \oplus H)$ is a subvariety of $X$ corresponding to a line
bundle
numerically equivalent to $T-wf$ while $\Bbb{P}(W)$ is an elliptic
curve $\Gamma$ on $X$, isomorphic to $C$, which is numerically equivalent
to
$T^2+xT f$ , for some $x \in \Bbb{Z}$. As the cycles
corresponding to
$\Bbb{P}(W)$ and $\Bbb{P}(G \oplus H)$ do not intersect, from $(T^2+xTf)
(T-wf)=0$ we get $x=-(g+h)$. If $D$ is a very ample line bundle on $X$,
$D_{|\Gamma}$ is very ample, hence $D  \Gamma = b+a\mu^-(E) \geq 3.$

ii) $E=H \oplus G$ where  rk$H=1$, rk$G=2$, $h=$ deg$H$,
$g= $ deg$G$.
As in  i), we get that $Z= \Bbb{P}(H)$ is numerically equivalent to
$T^2-gTf$
and the very ampleness of $D_{|Z}$ implies $b+ah \geq 3$. If $h \leq
\frac{g}{2}$ this
concludes the proof.
Otherwise let us denote by $Y$ the smooth surface $\Bbb{P}(G)$. As usual
Num($Y$)
is generated by $T'=T_{|Y}$ and $f'=f_{|Y}$. The very ampleness of $D$
implies the one
of
$D_{|Y}$ and by section 3, $D_{|Y}$  is very ample if and only if
\begin{itemize}
\item[]$b+a \frac {g}{2} \geq 3$ if $g$ is even and
\item[]$b+a \frac {g}{2} > 1$ if $g$ is odd.
\end{itemize}
If $g$ is even, a necessary condition for the very ampleness of $D$ is
$b+a\frac{g}{2} \geq 3$ i.e. $b+a\mu^-(E) \geq 3$  which is the desired
condition.

If $g$ is odd,  necessary conditions for the very ampleness of $D$ are
both $b+ah
\geq 3$ and $b+a\frac {g}{2} > 1.$ Hence  only the sufficient condition
$b+a\mu^-(E) \geq 3$ is obtained in this case.
\end{pf}
\section {rank $r$}

To deal with the case of $E$ vector bundle on an elliptic curve $C$, of
rank $r>3$  we
need to recall first the following
\begin{Thm}[\cite{H2}]
Let $E$ be a vector bundle of rank $r$ on an elliptic curve $C.$  $E$ is
ample if and
only if every indecomposable direct summand $E_i$ of $E$ has deg$E_i >0.$
\end{Thm}
Our method based on Lemma \ref{reductiontoS}, \ref {T+F} and \ref{T}
forces us to investigate first the case deg$E =3$, then deg$E=1, 2$ and
 finally deg$E\geq 4$.

\subsection{ deg$E=3$}
\begin{Thm}
\label{r4d3teo}
Let $E$ be a vector bundle over an elliptic curve $C$, with deg$E=3$ and
rank$E=4$,
and let $D\equiv aT+bf$ be a line bundle on $X=\Bbb{P}(E)$. Then the
following
conditions hold:
\begin{enumerate}
\item[i)] If $E$ is indecomposable, and $b+\frac{3}{4} a >\frac{3}{4}$,
$D$ is very ample
if and only if
$b+a \geq 3$
\item[ii)] If $E$ is decomposable and ample, and $b+ \frac{a}{3}
>\frac{1}{3}$  $D$ is
very ample if
$b+
\frac{a}{2} >2$.
 \end{enumerate}
\end{Thm}
\begin{pf}
i) If deg$E=h^0(\cal{T})=3$ then we can apply Lemma 4.3 and 4.1 with $A=
\cal T$,
 as
 $b+(a-1)\mu^-(E) = b+\frac{3}{4} a- \frac{3}{4} >0$ by
hypothesis. Hence $D$ is very ample if and only if $D_{|S} $ is very
ample, where $S$
is a suitable element of $|\cal{T}|$. There exists a vector bundle $E'$
with
rk$E'=3$, deg$E'=3$ , given by
\begin{equation}
\label{succdelT}
0\to \cal{O}_C \to E \to E' \to 0 ,
\end{equation}
such that  $S=\Bbb{P}(E')$ and Num(S)
is generated by
$T'$ and
$f'$, where
$T'= T_{|S}$ and $f'= f_{|S}$ so that $D_{|S}\equiv aT'+bf'$.

By section 4, $D_{|S}$ is very ample if and only if $b+a \geq 3$. Indeed
if $E'$ is
indecomposable the necessary and sufficient condition for the very
ampleness is $b+a
\geq 3$.
If $E'=\bigoplus_i E'_i$ then $\mu(E'_i) \geq
\frac{3}{4}$  for all $i.$ Hence the only possibilities for a
decomposition of $E'$ are:
 \begin{enumerate}
\item[1)] $E'=F \oplus G \oplus H$, with $F,G,H$ line bundles all of
degree $1.$
\item[2)] $E'= H \oplus G$, with rank $G=2$, rank $H=1$, deg $G=2$, deg
$H=1.$,
\end{enumerate}
Again by section 4 in both the  above cases, the necessary and sufficient
condition
for the very ampleness of $D_{|S}$ is $b+a \geq 3$.

ii) Let us suppose that $E$ is decomposable and ample. Then the possible
decompositions for $E$ give $\mu^-(E) =
\frac{2}{3},\frac{1}{3},\frac{1}{2}.$
Note that the condition $b+ \frac{a}{3} > \frac{1}{3}$ allows us to apply
lemma 4.3 with
$A=
\cal T$ in any case.  If $S$ is the usual element of $|\cal{T}|$,we get
that $D$ is very ample on $X$ if and only if $D_{|S}$ is. If
$S=\Bbb{P}(E')$ , $E'$ could
be decomposable. In this case  $\mu^-(E') \ge
\frac{1}{2}.$ Therefore the condition $b+ \frac{a}{2} > 2$ guarantees the
very
ampleness of $D_{|S}$  by section \ref{rango3sec}.
 \end{pf}
Because  $a\geq 2,$  Theorem \ref{r4d3teo} immediately gives the
following:
\begin{Cor}
\label{corond3r4}
Let $E$ be an  ample vector bundle over an elliptic curve $C$,
with deg$E=3$ and  rank$E=4$, and let $D \equiv aT+bf$ be a line bundle on
$X=\Bbb{P}(E)$ with $b+ \frac{a}{3} > \frac{1}{3}$. If  $b+ \frac{a}{2}
>2$ then  $D$ is
very ample.
\end{Cor}

\begin{Thm}
Let $E$ be an ample  vector bundle over an elliptic curve $C$, with deg
$E=3$ and
$r=$ rank$E\geq 4$, and let $D\equiv aT+bf$ be a line bundle on
$X=\Bbb{P}(E)$ with
$b+a\mu^-(E) >\frac{3}{5}$. Then $D$
is very ample if
$b+\frac{a}{2} >2$ and $b+\frac{a}{3} > \frac{1}{3}$.
\label{deg3anyr}
\end{Thm}
\begin{pf}
Proceed  by induction on $r=$ rank $E$. If $r=4$ the inductive
hypothesis is verified by Corollary \ref{corond3r4}. Let  $r \geq 5.$ It
is
$h^0(X, \cal T)=h^0(C, E)=3.$ Notice that $\mu^-(E) \le \mu(E)
\le\frac{3}{r}\leq\frac{3}{5}.$
Therefore  $b + (a-1) \mu^-(E) > 0$ and Lemma \ref{T} can be
applied, with $A = \cal{T}.$  Let
$S=\Bbb{P}(E')$ be as in (\ref{succdelT}) with deg $(E')=3$ and
rk $(E')=r-1\geq 4$.  Num$(S)$ is generated by $ T'=T_{|S}$ and $
f'=f_{|S},$ so that $
D_{|S}
\equiv aT'+bf'$. Notice that
$E'$ is ample being a quotient of $E.$
Notice that $\mu^-(E) \leq \mu^-(E')$. Indeed there exists a
map from at least one direct summand  $E_k$ of $E$ and  the summand $E'_j$
of $E'$
which realizes $\mu^-(E')$ and so we get $\mu^-(E) \leq \mu(E_k) \leq
\mu(E_j) =
\mu^-(E')$.
 As $b+a\mu^-(E') \geq b+a\mu^-(E) >\frac{3}{5}$, by induction
$D_{|S} $ is very ample if $b+\frac{a}{2} >2$,
$b+\frac{a}{3}>\frac{1}{3}$.
\end{pf}
\subsection{deg $E =2$}
\begin{Thm}
\label{d2anyrindecteo}
Let $E$ be an indecomposable vector bundle over an elliptic curve $C$,
with
deg$E=2$, rk $E= r \ge4.$  Let $D\equiv aT+bf$ be a line bundle on
$X=\Bbb{P}(E).$
 Then the following conditions hold:
\begin{enumerate}
\item[i)] If $r=4$ , $D$ is very ample if $b+\frac{a}{2}  >2$
\item[ii)] If  $r \ge 5 $, $D$ is very ample if $b+\frac{a}{2} >2$ ,
$b+\frac{a}{3}
>\frac{1}{3}$ and $b+a\frac{2}{r} > 1+ \frac{1}{r}$.
 \end{enumerate}
\end{Thm}
\begin{pf}
Let $A$ be as in Lemma \ref{T+F} and notice that $h^0(A) = 2 + r > $ deg
$E + 3.$ By
Remark \ref{criteriodelbutler} and the assumptions in i) and ii) we can
assume
\begin{equation}
\label{rangeconepsilon}
 b + a\frac{2}{r} = 1+\frac{\varepsilon}{r},
\end{equation} where
$\varepsilon = 1,2....r.$ If
$b+a\frac{2}{r} >1+\frac{2}{r}$ the assumptions of  Lemma
\ref{reductiontoS} and
\ref{T+F} are satisfied. If
 $b+a\frac{2}{r} = 1+\frac{2}{r}$ the line bundle $D-A$ has degree $0$ and
Lemma \ref{AA} shows that  Lemma \ref{reductiontoS} can be used if
\begin{equation}
\label{1+2sur}
b+a\frac{2}{r}
\geq 1+\frac{2}{r}.
\end{equation}

 Therefore condition (\ref{1+2sur}) can be rewritten using
(\ref{rangeconepsilon}) as
 \begin{equation}
\label{1+1sur}
b+a\frac{2}{r} > 1+\frac{1}{r}.
\end{equation}
 Let $S=\Bbb{P}(E')$, where $E'$ is as in (\ref{secondasucc}), where
 deg $ E'=3$ , rk $ E'=r-1\geq 3$, Num$(S)$ is generated by $ T'=T_{|S}$
and $
f'=f_{|S},$ and $D_{|S} \equiv aT'+bf'$. Notice that  $E'$ is
indecomposable or
decomposable and ample because
$ \mu^-(E') \geq \mu^-(E)=\frac{2}{r} > 0$.

To prove i) assume $r=4.$  In this case rk  $E'=3$ and deg $E'=3$. If $E'$
is
indecomposable, as $b+a \mu^-(E') >b+a\mu^-(E)=b+\frac{a}{2}> 1,$ by
Theorem
\ref{Thmrango3}
$D_{|S}$ is very ample if
$b+a
\geq 3$ .  By the same theorem if $E'$ is  decomposable and ample then
$D_{|S}$ is
very
 ample if $b+a \mu^-(E') > 2.$
The possible values for $\mu^-(E')$  are $1$ or $\frac{1}{2}.$ Therefore
 by Theorem \ref{Thmrango3} it follows that  $D_{|S}$ is very
ample  if
\begin{equation}
b+\frac{a}{2} >2.
\label{b+asu2}
\end{equation}   Putting (\ref{1+1sur}) and (\ref{b+asu2}) together it
follows that
in the case $r=4$
$D$ is very ample if
$b+\frac{a}{2} > \frac {5}{4}$ and  $b+\frac{a}{2} >2$ i.e. $b+\frac{a}{2}
>2$.

To prove ii) assume $r-1\geq 4.$ Since
 $b+a\mu^-(E')  \geq b+a\mu^-(E) =
b+a\frac{2}{r} >1+\frac{1}{r}>\frac{3}{5}$ from Theorem \ref{deg3anyr} it
follows that
$D_{|S}$ is very ample if
 $b+\frac{a}{2} >2$ and  $b+\frac{a}{3} >\frac{1}{3}.$ Therefore if $r\geq
5$, $D$ is very ample if $b+\frac{a}{2} >2$  ,  $b+\frac{a}{3}
>\frac{1}{3}$ and
 $b+a\frac{2}{r} > 1+ \frac{1}{r}$.
\end{pf}

\begin{Thm}
\label{d2anyrdecteo}
Let $E$ be a decomposable and ample vector bundle of rank $ r\ge 4$ over
an elliptic
curve
$C$, with deg $E=2.$
 Let $D\equiv aT+bf$ be a line bundle on
$X=\Bbb{P}(E).$ Then the following conditions hold:
\begin{enumerate}
\item[i)] If $r =4$ , $D$ is very ample if $b+\frac{a}{2} >2$  and
$b+a\mu^-(E) >
\frac{3}{2}$
\item[ii)] If $r \geq 5$, $D$ is very ample if $b+a\mu^-(E) >
1+\frac{2}{r} $,
$b+\frac{a}{3}>\frac{1}{3}$ and $b+ \frac{a}{2} > 2$.
\end{enumerate}
\end{Thm}
\begin{pf}
Let $A$ be as in Lemma \ref{T+F} and notice that $h^0(A) = r+2
\geq$ deg $E + 3.$ Also notice  that $\mu^-(E) \leq \mu(E) =
\frac{2}{r}$.  In our hypothesis we have  $b+(a-1)\mu^-(E)-1\geq
b+a\mu^-(E) -\frac{2}{r}-1 > 0$.

 Therefore Lemma \ref{reductiontoS} and \ref{T+F} can be applied.
Let $S=\Bbb{P}(E')$, where $E'$ is as in (\ref{secondasucc})
where deg $E'=3$ and $rk(E')=r-1 \geq 3$. If $r \geq 5$ we can apply
Theorem
\ref{deg3anyr} to
 $D_{|S}. $  Indeed $b+a\mu^-(E')  \geq b+a\mu^-(E)  > 1+ \frac{2}{r} >
\frac{3}{5}$.
Hence $D_{|S}$ is very ample if $b+\frac{a}{3}>\frac{1}{3}$ and $b+
\frac{a}{2} > 2.$

If $r =4$ we can apply Theorem \ref{Thmrango3}. If $E'$
is indecomposable then
$D_{|S}$ is very ample if
$b+a \geq 3.$ If $E'$ is decomposable, noticing that $E'$ is  ample, the
condition is $
b+a\mu^-(E')  > 2.$ In this case $\mu^-(E')$ can be $1$ or $\frac {1}{2}$
and hence the
worst sufficient condition becomes
$b+ \frac{a}{2} > 2$.
\end{pf}
 Theorem \ref{d2anyrindecteo} and \ref{d2anyrdecteo} give  the following
Corollary.
\begin{Cor}
\label{corond2}
Let $E$ be an ample  vector bundle over an elliptic curve $C$, with deg
$E=2$
 and let $D\equiv aT+bf$ be a line bundle on
$X=\Bbb{P}(E).$  Then the following conditions hold:
\begin{enumerate}
\item[i)] If  rank $E \ge 5$ , then $D$ is very ample
 if  $b+a\mu^-(E) >1+\frac{2}{r}$, $b+ \frac{a}{2} > 2$ and $b+
\frac{a}{3} > \frac{1}{3}.$
\item[ii)] If  rank $E=4$, $D$ is very ample if
$b+a\mu^-(E) > \frac{3}{2}$ and $b+ \frac{a}{2} > 2$.
\end{enumerate}
\end{Cor}


 \subsection{deg $E =1$}

Note that if deg$E=1$,  $E$ is ample if and only if it is  indecomposable.
\begin{Thm}
Let $E$ be an indecomposable vector bundle over an elliptic curve $C$,
with deg
$E=1$, rk $E = r  \geq 4$,  and let $D\equiv aT+bf$ be a line bundle on
$X=\Bbb{P}(E)$ with $b+ \frac{a}{r} >1$. Then the following conditions
hold:
\begin{enumerate}
\item[i)] If $r= 4,$  $D$ is very ample if  $b+\frac{a}{2} > 2$
\item[ii)] If $r = 5,$ $D$ is very ample if  $b+\frac{a}{3}>\frac{3}{2}$
and
$b+\frac{a}{2} > 2$.
\item[iii)] If $r \geq 6$, $D$ is very
ample if $b+\frac{a}{r-2} > 1+ \frac{2}{r - 1}$ and $b+\frac{a}{2} >2$.
\end{enumerate}
\end{Thm}
\begin{pf}
 In our hypothesis it is  $(a-1)\frac{1}{r} +b-1 > 0.$  Lemma \ref{AA}
allows us to
apply  Lemma \ref{reductiontoS} and \ref{T+F} when $(a-1)\frac{1}{r} +b-1
\ge
0,$ noticing that if  $A =
\cal{T}+\cal{F_P}$  it is
$h^0(A)=1+r
\geq 5$.  If $S=\Bbb{P}(E')$, where $E'$ is as in (\ref{secondasucc}),
notice that $E'$ is ample because $E$ is. It is
rk $E'=r-1\geq 3$ and
deg $E'=2.$ Let us now distinguish the three cases according to the values
of $r.$
Theorem \ref{deg3anyr} and Corollary \ref{corond2} will be used.
\begin{enumerate}
\item[i)] If rank$E=4$, rank$E'=3$,  and $D_{|S}$ is very ample if
$b+\frac{2}{3}a >\frac{4}{3} $ when $E'$ is
 indecomposable  while if $E'$ is  decomposable and ample  the condition
is
$b+a\mu^-(E') >2$ In the worst case $ \mu^-(E')= \frac{1}{2}$ so we have
$b+\frac{a}{2} >2$.
\item[ii)] If rank $E=5$ , rank $E'=4$, and $D_{|S}$ is very ample if
$b+a\mu^-(E') >
\frac{3}{2}$, and
 $b+\frac{a}{2} >2$.
 As $E'$ is indecomposable or decomposable and ample,
 in the worst case $ \mu^-(E')= \frac{1}{3}$ and so it is enough to ask
$b+\frac{a}{3} > \frac{3}{2}$ and $b + \frac{a}{2} > 2.$
 \item[iii)] If rank $E \geq 6$ then rank $E' \geq 5$, hence by Corollary
\ref{corond2}
 $D_{|S}$ is very ample if $b+a\mu^-(E') > 1+\frac{2}{r-1} $,
$b+\frac{a}{3}>\frac{1}{3}$
and
 $b+ \frac{a}{2} > 2$.
 If $E'$ is indecomposable then $\mu^-(E') = \frac{2}{r-1}$ while if
$E'$ is decomposable and ample, then  $\mu^-(E') = min (\frac{1}{s},
\frac{1}{r-1-s})$
with $ s= 1..r-2$. So the condition  $b+a\mu^-(E') > 1+\frac{2}{r-1} $ can
be
substituted by $b+\frac {a}{r-2} >1+\frac{2}{r-1} $   which implies the
condition $b+\frac{a}{3}>\frac{1}{3}.$
  \end{enumerate}
\end{pf}
%
%
\subsection{deg $E \geq 4$}

\begin{Thm}
Let $E$ be an ample  vector bundle over an elliptic curve $C$, with
deg$E=d \geq 4$,
rk $E = r \geq 4$ and $d < r.$  Let
$D\equiv aT+bf$ be a line bundle on $X=\Bbb{P}(E).$ If $b+
\frac{a}{d-1} > 2$ and $b + (a-1) \mu^-(E) > 0$ then
$D$ is very ample.
\end{Thm}
\begin{pf}
The proof proceeds by induction on $r.$
If $r=4$ the smallest possible value of $\mu^-(E)$ is $\frac{1}{3}$.
Since $d\ge 4$ it is
$\frac{1}{d-1} \le \frac{1}{3}.$ Therefore $b + a\mu^-(E) \ge b +
\frac{a}{3} \ge b +
\frac{a}{d-1} > 2.$  Lemma \ref{lemmadiEnrique} and Theorem
\ref{butlerteo}
conclude the proof of the initial inductive step. Let us suppose the
statement
true for $r-1$ and prove it for
$r$. Let $A$ be as in Lemma \ref{T}  and notice that $h^0(A) = d  \geq 4.$
Therefore
 Lemma
\ref{T} and
\ref{reductiontoS} can be applied.  By considering $S=\Bbb{P}(E')$ where
$E'$ is as in
(\ref{succdelT}),  we get that
$E'$ is ample with deg $E'=d \ge 4$,  rk $E'= r' =r-1 \ge 4.$ Because
$\mu^-(E') \ge \mu^-(E)$
it is
$b+ (a-1) \mu^-(E') \ge b + (a-1)\mu^-(E) >0$ and again $b + \frac{a}{d-1}
> 2.$
Hence we can conclude by induction hypothesis.
\end{pf}
\begin{rem}
Note that in the above theorem, if $E$ is indecomposable, by normalizing
$E$ we can always assume
$d < r.$
\end{rem}
In a very particular case we can say a little more:
\begin{Prop}
Let $E$ be an indecomposable vector bundle over an elliptic curve $C$,
 with deg $E=d \geq 4$, rk $E = d+1$,
and let $D\equiv aT+bf$ be a line bundle on $X=\Bbb{P}(E).$
If $b+(a-1)\frac{d}{d+1} >0$
and  $b+ a > 2$ then $D$ is very ample.
\end{Prop}
\begin{pf}
It is  $h^0(X, \cal{T}) = d \ge 4.$ Our hypothesis $b+(a-1)\frac{d}{d+1}
>0$
shows that Lemma
\ref{reductiontoS} and
\ref{T}  can be applied.

By considering  $S=\Bbb{P}(E')$ , where $E'$ is as in (\ref{succdelT}) it
is  deg $E'$ = rk
$E'$= $d$. A sufficient condition for the very ampleness of $D_{|S}$ is
$b+ a \mu^-(E') > 2$ by Lemma \ref{lemmadiEnrique} and Theorem
\ref{butlerteo}.
We claim that in this case
$\mu^-(E')=1$. If $E'$ is indecomposable it is
$\mu^-(E') = \mu(E') =1$. If $E'$
is decomposable we can  suppose that $E'= \oplus G_j$, with $r_j =  $ rk
$G_j
\geq 1$ , $d_j=$ deg $G_j \geq 1$ (as $\mu(G_j) \geq \mu(E)=\frac
{d}{d+1}$) and
$\sum d_j =\sum r_j =d$. Hence we have $\frac{d_j}{r_j} \geq \frac
{d}{d+1},$  for all
$j,$ which implies $r_j-d_j \leq \frac{d_j}{d} < 1,$  for all $ j.$ Hence
$r_j - d_j \leq 0 ,$
for all  $j$ i.e. $d_j=r_j+s_j$ with $s_j \geq 0,$  for all  $j.$ But $d=
\sum d_j = \sum
(d_j+s_j )= d+
\sum s_j$ and $\sum s_j =0$ i.e. $s_j=0,$  for all $ j$. Hence $\mu^-(E')
=1$.
\end{pf}

\end{document}